\documentclass{article}
\usepackage{amsmath}
\usepackage{theorem}
\usepackage{amssymb}
\usepackage{graphicx}
\usepackage{varioref}
\usepackage{subfigure}

\newtheorem{theorem}{Theorem}[section] 

\newtheorem{proposition}[theorem]{Proposition}
\newtheorem{corollary}[theorem]{Corollary}

{\theorembodyfont{\rmfamily}
\newtheorem{definition}[theorem]{Definition}
\newtheorem{example}[theorem]{Example} 
 
\newtheorem{remark}[theorem]{Remark}
\newtheorem{ack}{Acknowledgements}
}

\renewcommand{\bar}{\overline}
\newcommand{\A}{\ensuremath{\mathcal A}}

\newcommand{\D}{\ensuremath{\mathcal D}}

\newcommand{\X}{\ensuremath{\mathcal X}}
\newcommand{\C}{\ensuremath{\mathbb C}}
\renewcommand{\L}{\ensuremath{\mathsf{L}}}

\newcommand{\G}{\ensuremath{\mathcal G}}
\newcommand{\sG}{\ensuremath{\mathfrak{G}}}
\newcommand{\sK}{\ensuremath{\mathfrak{K}}}
\renewcommand{\NG}{\ensuremath{\mathcal {N}}}
\newcommand{\pV}{\ensuremath{\overline{V}}}
\newcommand{\R}{\ensuremath{\mathcal R}}
\newcommand{\V}{\ensuremath{\Sigma}}
\renewcommand{\Re}{\ensuremath{\mathbb R}}
\newcommand{\bR}{\ensuremath{\bar{\mathcal R}}}
\newcommand{\kk}{\ensuremath{\Bbbk}} 
\newcommand{\K}{\ensuremath{\Bbbk}}
 
\newcommand{\Cir}{\ensuremath{\mathcal C}}
\newcommand{\Z}{\ensuremath{\mathbb Z}}

\renewcommand{\P}{\ensuremath{\mathbb P}}
\newcommand{\bl}{\ensuremath{{\boldsymbol{\lambda}}}}
\renewcommand{\l}{\ensuremath{\lambda}} 

\newcommand{\be}{\ensuremath{{\boldsymbol{\eta}}}}
\newcommand{\e}{\ensuremath{\eta}}

\newcommand{\bx}{\ensuremath{{\boldsymbol{\xi}}}}
\newcommand{\x}{\ensuremath{\xi}}
\newcommand{\bm}{\ensuremath{{\boldsymbol{\mu}}}}
\newcommand{\m}{\ensuremath{\mu}}
\newcommand{\bn}{\ensuremath{{\boldsymbol{\nu}}}}
\newcommand{\n}{\ensuremath{\nu}}
\newcommand{\codim}{\operatorname{codim}}
\newcommand{\charac}{\operatorname{char}}
\newcommand{\supp}{\operatorname{supp}}

\newcommand{\depth}{\operatorname{depth}}

\newcommand{\qed}{\hfill \mbox{$\Box$}}
\newenvironment{proof}{\medskip\noindent {\it Proof:}}{\hfill \qed \medskip \break}

\begin{document}

\title{Resonance varieties over fields of \\ positive characteristic}
\author{Michael J. Falk}

\maketitle

\begin{abstract}
Let \A\ be a hyperplane arrangement, and \kk\ a field of arbitrary characteristic. We show that the projective degree-one resonance variety $\R^1(\A,\kk)$ of \A\ over \kk\ is ruled by lines, and identify the underlying algebraic line complex $\L(\A,\kk)$ in the Grassmannian $\G(2,\kk^n), n=|\A|.$  $\L(\A,\kk)$ is a union of linear line complexes corresponding to the neighborly partitions of subarrangements of \A. Each linear line complex is the intersection of a family of special Schubert varieties corresponding to a subspace arrangement determined by the partition.

In case \kk\ has characteristic zero, the resulting ruled varieties are linear and pairwise disjoint, by results of A.~Libgober and S.~Yuzvinsky. We give examples to show that each of these properties fails in positive characteristic. The (4,3)-net structure on the Hessian arrangement gives rise to a nonlinear component in $\R^1(\A,\overline{\Z}_3),$ a cubic hypersurface in $\P^4$ with interesting line structure. This provides a negative answer to a question of A.~Suciu. The deleted $B_3$ arrangement has linear resonance components over $\Z_2$ that intersect nontrivially. 

\end{abstract}

\begin{section}{Resonance and characteristic varieties}

Arising out of the study of local system cohomology and fundamental groups, characteristic and resonance varieties of complex hyperplane arrangements have become the object of much of the current research in the field. The study of resonance varieties in particular has led to surprising connections with other areas of mathematics: generalized Cartan matrices \cite{LY00}, Latin squares and loops, nets \cite{Yuz04}, special pencils of plane curves \cite{LY00,FY06}, graded resolutions and the BGG correspondence \cite{EPY,SchSuc,SchSuc2}, and the Bethe Ansatz \cite{CV}. In this paper we establish a connection between resonance varieties and projective line complexes, which becomes evident only when one works over a field of positive characteristic.

Let \A\ be a (central) complex hyperplane arrangement of rank $\ell,$ \kk\ a field, and $A=A_\kk(\A)$ the Orlik-Solomon algebra of \A\ over \kk. Since $A$ is a quotient of an exterior algebra, left-multiplication by a fixed element $a\in A^1$ defines a chain complex $(A,a)$:$$0\to A^0\to A^1\to \cdots \to A^\ell\to 0.$$
The $d^{\sl{th}}$ {\em resonance variety} of $\A$ over \kk\ is $$\R^d(\A,\kk)=\{a\in A^1 \ |  H^d(A,a)\not =0\}.$$ $\R^d(\A,\kk)$ is filtered by $\dim_\kk(H^d(A,a)).$
Let $\R(\A,\kk)=\R^1(\A,\kk).$ An element $a\in A^1$ lies in $\R(\A,\kk)$ if and only if there exists $b\in A^1$ not proportional to $a$ with $ab=0.$ Identifying $A^1$ with $\kk^n$ using the canonical basis corresponding to the hyperplanes of \A, $\R(\A,\kk)$ is seen to be a homogeneous affine algebraic variety. Hence $\R(A,\kk)$ determines a projective variety $\bR(A,\kk)$ in $\P^{n-1},$ which we call the (degree-one) projective resonance variety of \A\ over \kk.

If \kk\ has characteristic zero, $\bR(A,\kk)$ has remarkable algebraic and combinatorial features. It consists of disjoint linear subspaces, each of which corresponds to a combinatorial structure, called  a {\em multinet}, on a subarrangement of \A. If $\kk=\C,$ the multinet structure arises from a pencil of plane projective curves whose singular elements include at least three completely reducible, not necessarily reduced curves \cite{FY06}, see also \cite{Di06,Lib06}. Furthermore, by \cite{LY00}, the stratification of $\R(\A,\kk)$ coincides with the stratification by dimension; that is,  the dimension of $H^1(A,a))$ is equal to the dimension of the irreducible component of $\bR(A,\kk)$ containing $a.$ 

The present work generalizes the description of $\bR(\A,\kk)$ from \cite{LY00} to arbitrary fields. In general, $\bR(\A,\kk)$ is a union of ruled varieties determined by neighborly partitions of \A. The partitions need not correspond to multinets, and the corresponding varieties need not be linear or disjoint, as our Examples~\ref{hessian} and \ref{deleted1} show. Example \ref{nonfano}, due to D.~Matei and A.~Suciu \cite{MatSuc1}, shows that the stratification of $\bR(\A,\kk)$ need not coincide with stratification by dimension.

The ruled variety corresponding to a partition $\Gamma$ of \A\ is easily described in terms of Schubert varieties. Let $\bar{K}$ be the subspace of $\P^{n-1}$ given by the kernel of the incidence matrix corresponding to codimension-two flats meeting more than one block of $\Gamma.$ The subspaces of $\bar{K}$ given by the kernels of characteristic functions of the blocks are called {\em directrices}. Each of these subspaces determines a Schubert variety in the Grassmannian of lines in $\P^{n-1},$ consisting of the lines it meets. The intersection of these Schubert varieties - the set of lines that meet every directrix - is the line complex $\L(\Gamma,\kk)$ determined by $\Gamma,$ and the union $\pV(\Gamma,\kk)$ of the lines in $\L(\Gamma,\kk)$ is the constituent of $\bR(A,\kk)$ arising from $\Gamma.$ In every example $\pV(\Gamma,\kk)$ is an irreducible component of $\bR(A,\kk)$ but we have no proof of this in general. For this reason we will continue to refer to $\pV(\Gamma,\kk)$ as a ``constituent" of $\bR(A,\kk)$. By \cite{Yuz95}, $\pV(\Gamma,\kk)$ is contained in the diagonal hyperplane $\bar{K}_0$ in $\bar{K}.$

If \kk\ has characteristic zero, the results of \cite{LY00} imply that the directrices are hyperplanes in $\bar{K},$ and therefore $\pV(\Gamma,\kk)=\bar{K}.$ 
In order for $\pV(\Gamma,\kk)$ to be nonlinear, there must be at least three directrices which are not hyperplanes, none of which is a point. In particular, $\bar{K}_0$ must have dimension at least three. In characteristic zero, $\dim(\bar{K}_0)\geq 3$ implies the multinet has at least 4 blocks, and only one example (up to isomorphism) is known - the $(4,3)$-net structure \cite{Yuz04} on the Hessian arrangement of twelve lines determined by the nine inflection points of a smooth cubic in $\P^2(\C).$ In positive characteristic there are a few examples known for which $\dim(\bar{K}_0)=3,$ but they all have directrices which are points. For the Hessian configuration, and $\kk=\bar{\Z}_3,$ $\bar{K}_0$ is a $\P^4,$ and the four directrices are planes lying in special position. The ruled variety $\pV(\Gamma,\kk)$ is an irreducible, singular cubic threefold in $\P^4$ with very interesting geometry. This might be the only example of a nonlinear resonance component in nature.

All of this is depends only on the underlying matroid of \A, and can be carried out for any matroid, realizable or not. It would be interesting to examine resonance varieties associated with other finite affine or projective planes (the Hessian arrangement is a copy of $AG(2,3)$), which are not realizable over fields of characteristic zero.

In the rest of the introduction, we sketch some background and motivation for this research. The Orlik-Solomon algebra $A$ is isomorphic to the cohomology ring of the complement $X=\C^\ell -\bigcup \A$ with coefficients in \kk. For $\kk=\C$ the complex $(A,a)$ corresponds to a subcomplex of the twisted DeRham complex associated with a rank-one complex local system on $X$ corresponding to $a.$ The study of resonance varieties, over $\C$ or $\Z_p,$ is partly motivated by their connection with jumping loci for the cohomology of such local systems via the tangent cone theorem. 

For simplicity we will confine our discussion to the arrangement setting. Suppose \A\ consists of $n$ hyperplanes. 
Rank-one local systems on $X$ are parametrized by the complex torus $(\C^*)^n,$ with ${\bf t}=(t_1,\ldots, t_n)\in (\C^*)^n$ corresponding to the local system $\mathcal{L}_{\bf t}$ with monodromy around a loop about the $i^{\rm th}$ hyperplane given by multiplication by $t_i.$ The $d^{\rm th}$ {\em characteristic variety} of $X$ is $$\V^d(X,\K)=\{{\bf t}\in (\C^*)^n \ | \ \dim H^d(X,{\mathcal L}_{\bf t})\not =0\}.$$ Again we focus on the case $d=1.$

The stratification of $(\C^*)^n$ by $\dim(H^1(X,{\mathcal L}_{\bf t}))$ determines the first betti numbers of finite abelian covers of $X$, among which is  the Milnor fiber of the (non-isolated) singularity of $\bigcup\A$ at the origin.
Characteristic varieties over finite fields were first considered in \cite{Matei,MatSuc1,MatSuc2}. These determine further numerical invariants of $\pi=\pi_1(X)$, including $p$-torsion ($p\not = \charac(\K)$) in the first homology of finite-index subgroups. Resonance varieties over $\Z_p$ determine the same invariants for the nilpotent quotient $\pi/\pi^{(3)}.$ The cohomology of $(A,a)$ over $\Z_p$ gives upper bounds on the first betti numbers of cyclic covers of $X,$ such as the Milnor fiber \cite{CO98}.

The tangent cone theorem of Cohen and Suciu \cite{CS3} states that the degree-one resonance variety $\R(\A,\C)$ coincides with the tangent cone at the identity to $\V^1(X,\C).$ (There are now many generalizations of this theorem in different directions, resulting in applications of resonance varieties in several contexts - see e.g., \cite{CO00,L6,CS06a,DPS05}.) Thus the $\C$-resonance varieties determine the (subtorus) components of $\V^1(X,\C)$ passing through the identity, which therefore depend only on the underlying matroid. It is not known whether the same holds true for $\V^1(X,\C)$ itself. According to Arapura \cite{Arap}, components of $\V^1(X,\C)$ are positive dimensional tori translated by torsion elements, or isolated points on $(S^1)^n.$ Both phenonema can occur for hyperplane complements \cite{CS3,Suc,Suc2,Coh2}. There is evidence of a connection between these translated components and resonance over $\Z_p$ - see Examples~\ref{nonfano} and \ref{deleted1}. Indeed, resonance varieties over $\Z_p,$ considered as a subgroup of $\C^*,$ exponentiate to rational points in $\V^1(X,\C),$ by the result of \cite{CO98} cited above. 

The tangent cone result fails in positive characteristic  - see \cite[Example 10.7]{Suc2}. It may be that resonance varieties over $\Z_p$ provide an obstruction of $\Z_p$-formality. There is at least a connection between $\bR(\A,\Z_p)$ and Massey products over $\Z_p,$ see Example~\ref{deleted1} and \cite{Mat05}. 

One can also study resonance over $\Z$ or $\Z_N$ with $N$ composite, and there are again empirical connections with $\V^1(X,\C).$ Our approach via neighborly partitions can be generalized to allow coefficients in an arbitrary commutative ring - see \cite{Fa04}.

Here is an outline of this paper. In Section \ref{resonance} we describe our main objects of study and recall the description of $\R(\A,\kk)$ from \cite{F7}, a decomposition of $\R(\A,\kk)$ parametrized by neighborly partitions. 
In Section \ref{sect-lines} we describe the geometry of the constituents of $\R(\A,\kk)$  in terms of projective line geometry. We exhibit Schubert calculus formulae that may be used in some cases to compute the dimension and degree of these varieties.  In Section~\ref{sect-examples} we examine the resonance varieties of the non-Fano, the deleted $B_3,$ and the Hessian arrangement, over fields of characteristic 2 and 3. Here we observe the phenomena special to the positive characteristic case discussed earlier. We remark on the connections between modular resonant weights and torsion points in the characteristic variety.

This article is a revision of part of the unpublished preprint \cite{Fa04}.
\end{section}

\begin{section}{The combinatorial decomposition of $\R(\A,\kk)$}
\label{resonance}

Let $\A=\{H_1,\ldots,H_n\}$ be a central arrangement of $n$ distinct linear hyperplanes in $\C^\ell$. The combinatorial structure of \A\ is recorded in the underlying matroid $\sG=\sG(\A).$ This is the matroid on $[n]:=\{1,\ldots, n\}$  whose set of circuits $\Cir$ consists of the minimal sets $C\subseteq [n]$ satisfying $\codim (\bigcap_{i\in C} H_i)<|C|$. Note that $\sG$ has no circuits of size one or two. From the topological standpoint, we are mainly interested in the topology of the complement $X(\A)=\C^\ell - \bigcup_{i=1}^n H_i$, which is determined to a large, albeit ultimately unknown extent by the underlying matroid $\sG$.

Let $\K$  be a field. Let $E_\kk(n)$ denote the free graded exterior algebra over $\K$ generated by $1$ and degree-one elements $e_i$ for $i\in [n]$. 

\begin{definition}
The {\em Orlik-Solomon (OS) algebra} $A=A_\kk(\sG)$ of $\sG$ is the quotient of $E_\kk(n)$ by the homogeneous ideal $$I=(\partial e_C \ | \ C\in \Cir),$$ where $\partial$ is the usual boundary operator: $\partial e_C=\sum_{k=1}^p (-1)^{k-1}e_{i_1}\cdots\hat{e}_{i_k}\cdots  e_{i_p},$ for $C=\{i_1,\ldots,i_p\}$. 
\end{definition}
The image of $e_i$ in $A$ is denoted $a_i$. Then $A$ is a graded-commutative $\K$-algebra, generated by $1$ and the degree-one elements $a_i, 1\leq i\leq n.$ According to \cite{OT}, the dimension of $A^d_\kk(\sG)$ is independent of $\K$. In fact these dimensions are given by the coefficients of the characteristic polynomial of the lattice of flats of $\sG$. The Orlik-Solomon algebra is isomorphic to the cohomology algebra of the complement $X(\A)$ with coefficients in $\K$ \cite{OS1,OT}. The generators $a_i$ correspond to logarithmic 1-forms $d\alpha_i/\alpha_i$ where $\alpha_i: \C^\ell \to \C$ is a linear form with kernel $H_i$. 

\begin{example} Consider the arrangement of the six planes $x+ y=0,\  x-y=0,\  y+ z=0,\  y-z=0,\  z+x=0,$ and $z- x=0$ in $\C^3.$ This is the reflection arrangement of type $D_3,$ equivalent via a linear change of coordinates to the rank-three braid arrangement \cite{OT}. Labelling the elements of \A\ in the order given, the set \Cir\ of circuits of the underlying matroid is given by $\Cir=\{136, 145, 235, 246, 1234, 1256, 3456\}. $ These are the circuits of the cycle matroid $\sK_4$ of the complete graph on four vertices. The OS algebra $A$ has Hilbert series $H(A,t)=1 + 6t + 11t^3 + 6t^3.$ The circuit $136$ gives rise to the relation $\partial e_{136}=e_3 e_6-e_1 e_6 + e_1 e_3 =(e_1-e_3)(e_3-e_6)\in I,$ so that $(a_1-a_3)(a_3-a_6)=0$ in $A$.
\label{braid1}
\end{example}

\subsubsection*{Resonant weights} Let $\R(\A,\kk)=\{a\in A^1 \mid ab=0\ \text{for some} \ b\in A^1 -\kk a \},$ as in the introduction. We review the description of $\R(\A,\kk)$ from \cite{F7}. That paper treated the case $\kk=\C,$ but the methods apply to any coefficient field. The other approaches to resonance varieties in the literature (e.g., \cite{CS3,LY00,CO00}), apply only to the characteristic-zero case. 

For $\bx\in \kk^n$ set $a_\bx=\sum_{i=1}^n \x_i a_i$. Let $$Z(\bl,\K)=\{\be \in \kk^n \ | \ a_\bl  a_\be =0\}.$$ Usually $Z(\bl,\K)$ is abbreviated to  $Z(\bl)$, when no ambiguity results. If $\be\in Z(\bl)$ is not a scalar multiple of \bl, we call $(\bl,\be)$ a {\em resonant pair}. We say \bl\ is a {\em resonant weight} if \bl\ belongs to a resonant pair. The {\em support} $\supp(\bl,\be)$ of a resonant pair  is the set $\{i \in [n] \ | \ \l_i\not = 0 \ \text{or} \ \e_i\not = 0\}.$
Identifying $A^1$ with $\kk^n$ via the basis $\{a_1,\ldots, a_n\}$, $\R(\A,\kk)$ is precisely the set of resonant weights, i.e., $$\R(\A,\kk)=\{\bl\in \kk^n \mid \dim_\kk(Z(\bl))\geq 2\}.$$

The natural decreasing filtration of $\R(\A,\kk)$ is then defined by $$\R_k(\A,\K)=\{\bl \in \kk^n \ | \ \dim_\kk Z(\bl)\geq k+1\}, \ \text{for} \ k\geq 1.$$

\begin{example} Referring to Example \ref{braid1}, we have that $\bl=(1,0,-1,0,0,0)$ is a resonant weight for the braid arrangement, for any field, since $a_\bl a_\be=0$ for $\be=(0,0,1,0,0,-1)$, and \bl\ and \be\ are not parallel. In a similar way, every circuit of size three gives rise to a resonant pair.
\label{local1}
\end{example}
 
For $\bx=(\x_1,\ldots,\x_n) \in \kk^n$ and $S\subseteq [n]$, we denote the coefficient sum $\underset{i\in S}{\sum} \x_i$ by $\x_S$. We say \bl\ and \be\ are {\em parallel} if the matrix $\begin{bmatrix}\bl | \be \end{bmatrix}$ has rank less than two. Recall the {\em rank} of \A\ is the codimension of $\bigcap_{i=1}^n H_i$ in $\C^\ell.$

\begin{proposition}[\cite{Yuz95,F7}] Suppose \A\ has rank two, and $\bl,\be \in \kk^n$. 
\begin{enumerate}
\item If $n\geq 3,$ then $a_\bl a_\be=0$ 
if and only if $\l_{[n]}=\e_{[n]}=0.$
\item If $n=2,$ then $a_\bl  a_\be=0$ if and only if \bl\ and \be\ are parallel.
\end{enumerate}
\label{ranktwo}
\end{proposition}

Let \A\ be a central arrangement of arbitrary rank. Recall $\sG$ is the underlying matroid of \A, with ground set $[n].$ A {\em line} in $\sG$ is rank-two flat, or, equivalently, the closure in \sG\ of a two-element subset of $[n].$  Thus a line in $\sG$ corresponds to a maximal subarrangement of \A\ intersecting in a codimension-two subspace.
A line $X$ in $\sG$ is {\em trivial} if $|X|=2$. We denote the set of lines in $\sG$ by $\X=\X(\sG)$, and the set of nontrivial lines by $\X_0=\X_0(\sG)$. For $X \in \X$ let $A_X$ denote the subalgebra of $A$ generated by $\{a_i \mid i\in X\}$. Then $\A_X$ is isomorphic to $A(\A_X)$, where $\A_X=\{H\in \A \mid H\supseteq X\}.$ Moreover,  a fundamental result of Brieskorn affords us a direct sum decomposition $$A^2\cong \underset{X\in \X}{\oplus} A^2_X.$$  See \cite{OT} for proofs of these results. Proposition \ref{ranktwo} then yields a characterization of resonant weights for \A.
For $S\subseteq [n]$ define the {\em restriction} $\bx(S)$ of \bx\ to $S$ to be the element $(\x_i \ | \ i\in S)$ of $\K^{|S|}$.

\begin{theorem}[\cite{F7,LY00}] Suppose $\A$ is an arrangement of arbitrary rank. Then $\be\in Z(\bl)$ if and only if, for every $X\in \X(\sG)$, either 
\begin{enumerate}
\item $\bl(X)$ and $\be(X)$ are parallel, or
\item $X\in \X_0(\sG)$ and  $\l_X=\e_X=0.$ 
\end{enumerate}
\label{res-pair}
\end{theorem}

\begin{proof} Note that $a_\bl a_\be=\sum_{i<j} \begin{vmatrix}\l_i&\e_i\\ \l_j&\e_j\end{vmatrix}  \ a_i  a_j.$  By the direct sum decomposition above, $\be\in Z(\bl)$ if and only if $\be(X) \in Z(\bl(X))$ for every $X\in \X.$ The result then follows immediately from \ref{ranktwo}.  
\end{proof}

\subsubsection*{Neighborly partitions} Suppose $(\bl,\be)$ is a resonant pair. Define a graph $\Gamma=\Gamma_{(\bl,\be)}$ with vertex set  $[n]$, and $\{i,j\}$ an edge of $\Gamma $ if and only if  $\begin{vmatrix}\l_i&\e_i\\ \l_j&\e_j\end{vmatrix}=0$. If $\{i,j\}$ is an edge of $\Gamma$ we write $\{i,j\}\in\Gamma$. Note that $\{i,j\}\in \Gamma$ for every trivial line $\{i,j\}\in \X(\sG),$ by Prop.~\ref{ranktwo}(ii). Also, if $i\not \in \supp(\bl,\be)$ then $i$ is a cone vertex in $\Gamma$, that is, a vertex adjacent to every other vertex. A {\em clique} in $\Gamma$ is a set of vertices that are pairwise adjacent. Since \bl\ and \be\ are not parallel, $[n]$ itself is not a clique of $\Gamma$. By Theorem \ref{res-pair} and the definition of $\Gamma_{(\bl,\be)}$, we have the following.
 
\begin{corollary} If $X\in \X_0(\sG)$ is not  a clique of $\Gamma_{(\bl,\be)}$, then $\l_X=\e_X=0$.
If $X\in \X(\sG)$ is a clique in $\Gamma_{(\bl,\be)}$, then $\bl(X)$ is parallel to $\be(X)$.
\label{clique}
\end{corollary}

We define a {\em block} of $\Gamma$ to be a maximal clique. The blocks of $\Gamma$ cover $[n]$, but need not be disjoint. A cone vertex of $\Gamma$ is contained in every block.
\begin{definition} A graph $\Gamma$ with vertex set $[n]$ is {\em neighborly,} or more precisely {\em $\sG$-neighborly,} if for every $X\in \X(\sG)$ and every block $S$ of $\Gamma$, $|X\cap S|\geq |X|-1$ implies $X\subseteq S$.
\label{neighborlydef}
\end{definition}

Observe that a \sG-neighborly graph must include among its edges all the trivial lines $\{i,j\}\in \X(\sG)$. Also, if $i$ is a cone vertex of $\Gamma$, then $\Gamma$ is $\sG$-neighborly if and only if the induced subgraph on $[n]-\{i\}$ is $(\sG-i)$-neighborly.

\begin{theorem} Let $\Gamma=\Gamma_{(\bl,\be)}$ be the graph associated with a resonant pair of weights. Then $\Gamma$ is neighborly.
\label{neighbor}
\end{theorem}

\begin{proof}
Let $X\in \X(\sG),$ and $i\in X$ with $X-\{i\}$ a clique of $\Gamma.$  Suppose $X$ is not  a clique. Then $\l_X=\e_X=0,$ by Corollary~\ref{clique}. Let $j\in X$ with $\begin{vmatrix}\l_i&\e_i\\ \l_j&\e_j\end{vmatrix}\not=0.$ We have  $\begin{vmatrix}\l_X&\e_X\\ \l_j&\e_j\end{vmatrix}=0.$ But $\begin{vmatrix}\l_k&\e_k\\ \l_j&\e_j\end{vmatrix}=0$ for $k\in X-\{i\}$ since $j\not =i$ and $X-\{i\}$ is a clique. We conclude $\begin{vmatrix}\l_i&\e_i\\ \l_j&\e_j\end{vmatrix}=0,$ a contradiction. 
\end{proof}

The cone vertices of $\Gamma_{(\bl,\be)}$ are precisely the elements of $[n]-\supp(\bl,\be).$ Indeed, if $(\l_i,\e_i)\not=(0,0)$ is parallel to each $(\l_j,\e_j),$ then $\begin{bmatrix}\bl\ | \be \end{bmatrix}$ has rank one, so $\bl$ and $\be$ are parallel. By the same reasoning, the restriction of $\Gamma_{(\bl,\be)}$ to $\supp(\bl,\be)$ is transitive, i.e., every component is a clique. The resulting partition of $\supp(\bl,\be)$ is a neighborly partition in the sense of \cite{F7}. Our introduction of neighborly graphs on $[n]$ is merely a convenience to avoid the need to pass to subarrangements. They play a more crucial role when working over arbitrary commutative rings \cite{Fa04}.

By way of example, we illustrate in Figure~\ref{neighborlyfig} a pair of neighborly graphs with full support, depicted as partitions of the underlying matroid. Hyperplanes are represented by points; nontrivial rank-two flats by lines. Blocks of $\sG$-neighborly graphs are labelled by capital letters.

\begin{figure}[ht]
\centering
\subfigure[The cycle matroid $\sK_4$]{\includegraphics[height=4cm]{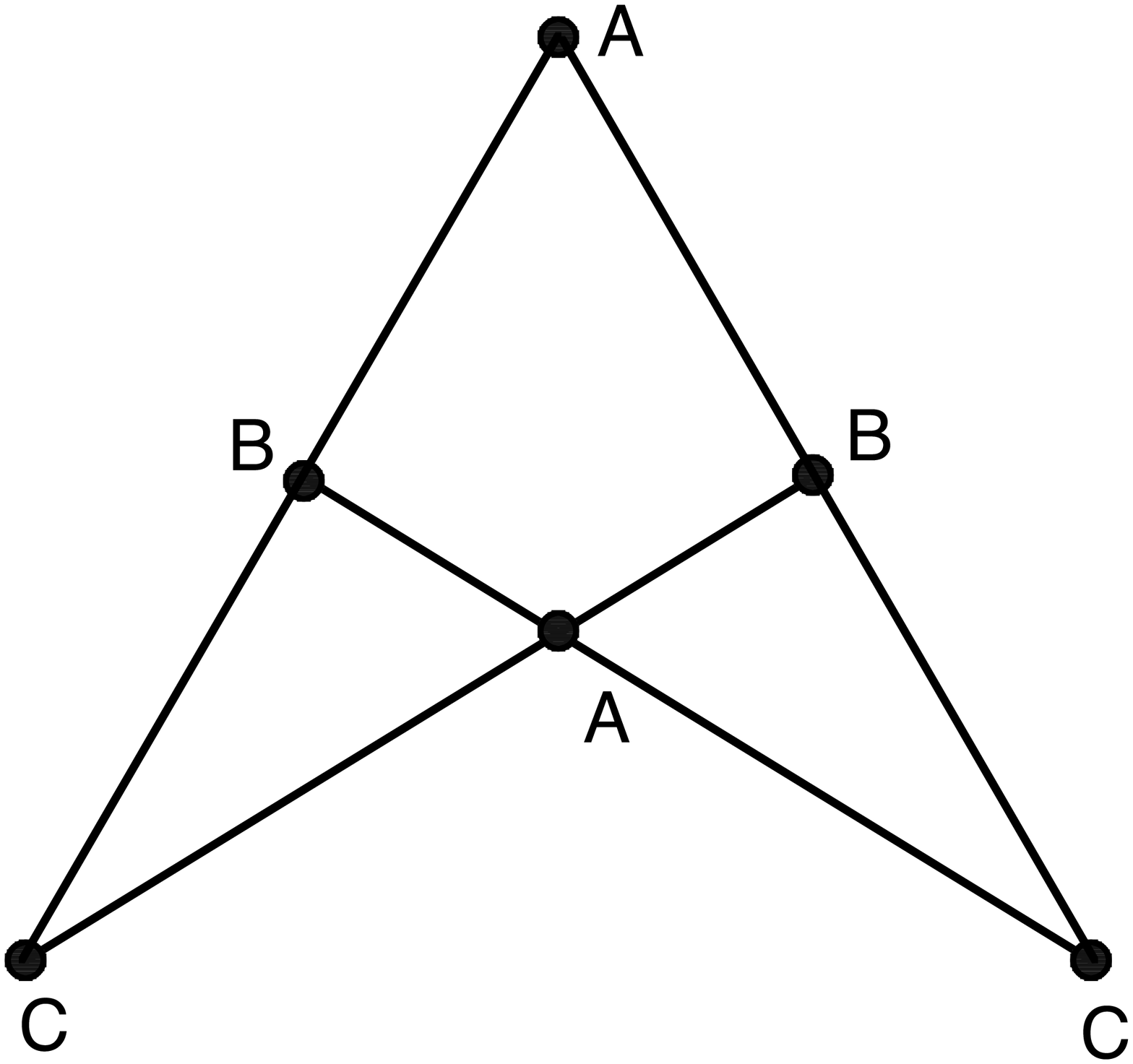}}
\subfigure[The non-Fano plane]{\includegraphics[height=4cm]{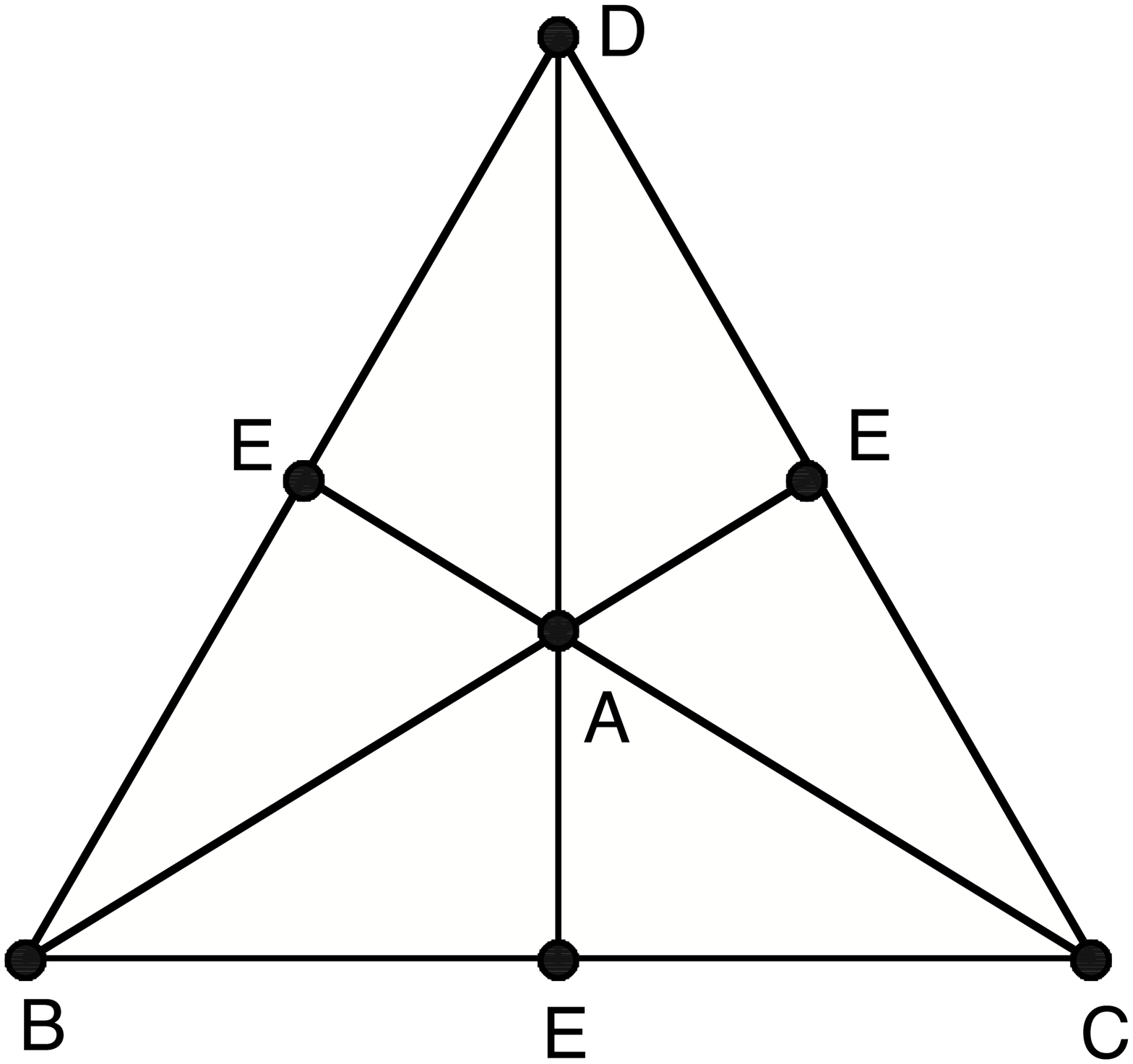}}
\caption{Neighborly partitions}
\label{neighborlyfig}
\end{figure}

\subsubsection*{Decomposition of $\R(\A,\kk)$}
For $\Gamma$ an arbitrary graph with vertex set $[n]$, set $$\X_\Gamma(\sG)=\{X\in \X_0(\sG) \ | \ X \ \text{is not a clique of} \ \Gamma\}$$ and 
$$K=K(\Gamma,\K)=\{\bx \in \kk^n \ | \ \x_X=0  \ \text{for every} \ X\in \X_\Gamma(\sG)\}.$$ 

For $\bl\in K$ we define  $$Z_\Gamma(\bl,\kk)=\{\be\in K \mid \begin{vmatrix}\l_i&\l_j\\ \e_i&\e_j\end{vmatrix} =0 \ \text{for every} \ \{i,j\}\in \Gamma\}.$$

In particular, if $\be\in Z_\Gamma(\bl,\K)$, and $S$ is a clique of $\Gamma$, then $\be(S)$ is parallel to $\bl(S)$. The converse may not be true, that is, the graph $\Gamma_{(\bl,\be)}$ may have more edges than the original graph $\Gamma$. We will write $Z_\Gamma(\bl,\K)$ as $Z_\Gamma(\bl)$ when it is not ambiguous. Note that $Z_\Gamma(\bl)$ is a linear subspace of $K.$

\begin{corollary} $Z_\Gamma(\bl)\subseteq Z(\bl)$.
\label{graph-cycles}
\end{corollary}

\begin{proof}
Let $\be\in Z_\Gamma(\bl).$ If $\bl(X)$ and $\be(X)$ are not parallel, then $X\in \X_\Gamma(\sG)$ and $\l_X=\e_X=0$. Then $\be\in Z(\bl)$ by Theorem \ref{res-pair}.
\end{proof}

With some care we can define a single graph $\Gamma$, depending only on $\bl$, such that $Z(\bl)=  Z_\Gamma(\bl).$ Let $E$ be a field extension of $\K$. An element $\bm \in E^n$ is called a {\em generic partner} of \bl\ if the following conditions are satisfied:

\begin{enumerate}
\item if $X\in \X_0(\sG)$ and $\m_X=0$, then $\e_X=0$ for all $\be \in Z(\bl,E)$, and
\item if $1\leq i<j\leq n$ and $\begin{vmatrix}\l_i&\l_j\\ \m_i & \m_j\end{vmatrix}=0$, then $\begin{vmatrix}\l_i&\l_j\\ \e_i & \e_j\end{vmatrix}=0$ for every $\be\in Z(\bl,E)$.
\end{enumerate}

Every $\bl\in \kk^n$ has a generic partner in $(\bar{\kk})^n,$ by Hilbert's Nullstellensatz. For computational purposes it may be desirable to find generic partners without passing to the algebraic closure of \kk. 

\begin{proposition}If $E$ is a field extension of $\K$ satisfying
$$|E|>\binom{n}{2}+|\X_0(\sG)|+1,$$ then every $\bl\in \R(\A,\kk)$ has a generic partner in $E^n$. 
\end{proposition}

\begin{proof}
Let $d=\dim_E(Z(\bl,E)).$ Then $d=\dim_\kk(Z(\bl,\kk)),$ and, in particular, $d\geq 2$ by assumption. The linear equations $\x_X=0$ determine $|\X_0(\sG)|$ hyperplanes in $Z(\bl,E)$. For $\bl$ fixed, the equations $\begin{vmatrix}\l_i&\l_j\\ \x_i & \x_j\end{vmatrix}=0$ also define hyperplanes in $Z(\bl,E)$, at most $n \choose 2$ of them. There are $|E|$ multiples of \bl\ in $Z(\bl,E).$ Since
\begin{eqnarray*}
|Z(\bl,E)|&=&|E|^d\\
&>&|E|^{d-1}\bigl({n\choose 2}+|\X_0(\sG)|+1\bigr)\\
&\geq& |E|^{d-1}\bigl({n\choose 2}+|\X_0(\sG)|)+|E|\\
&\geq& {n\choose 2}+|\X_0(\sG)|+|E|,
\end{eqnarray*}
there is a point $\bm$ of $Z(\bl,E)-E\bl$ missing the aforementioned subspaces. Such a \bm\ is a generic partner of \bl.
\end{proof}

We define $\Gamma_\bl=\Gamma_{(\bl,\bm)}$, where \bm\ is a generic partner of \bl\ over some extension field $E$ of \kk. It follows from condition (ii) that $\Gamma_\bl$ is well-defined. Note $\supp(\bl,\bm)\supseteq\supp(\bl,\be)$ for every $\be\in Z(\bl,\K)$. 

\begin{proposition} Suppose $\bl\in\R(\A,\kk).$ Then $Z(\bl)=Z_{\Gamma}(\bl)$ for $\Gamma=\Gamma_\bl.$ 
\label{goodgraph}
\end{proposition}

\begin{proof} We have $Z_\Gamma(\bl)\subseteq Z(\bl)$ by Corollary~\ref{graph-cycles}.  Let $\be\in Z(\bl),$ let $\bm$ be a generic partner of $\bl,$ so that $\Gamma=\Gamma_\bl=\Gamma_{(\bl,\bm)}.$ If $X\in \X_{\Gamma}(\sG),$ then $\bl(X)$ and $\bm(X)$ are not parallel, so $\l_X=0=\m_X$ by Theorem \ref{res-pair}. 
Then $\e_X=0$ by genericity of \bm. Hence $\be\in K(\Gamma,\K)$.
Also, by (ii) above,  $\begin{vmatrix}\l_i&\e_i\\ \l_j&\e_j\end{vmatrix}=0$ for every $\{i,j\}\in \Gamma$. Thus $\be\in Z_{\Gamma}(\bl)$. 
\end{proof}

\begin{corollary} Suppose $\bl\in\R(\A,\kk).$ Then $\dim(K(\Gamma_\bl,\kk))\geq 2.$
\label{dimen}
\end{corollary}

We can now decompose $\R(\A,\kk),$ and, more generally, $\R_k(\A,\kk)$ for each $k\geq 1.$ Let $\NG(\sG)$ denote the set of $\sG$-neighborly graphs with vertex set $[n]$, and $$\NG(\sG,\K)=\{\Gamma \in \NG(\sG) \ | \ \dim K(\Gamma,\kk)\geq 2\}.$$ 

For $\Gamma \in \NG(\sG,\K),$ set $$V_k(\Gamma,\K)=\{\bl \in K(\Gamma,\K) \ | \ \dim Z_\Gamma(\bl)\geq k+1\}.$$ 
Set $V_1(\Gamma,\K)=V(\Gamma,\K)$.

\begin{theorem}  $$\R_k(\A,\K)=\bigcup_{\Gamma\in \NG(\sG,\K)} V_k(\Gamma,\K).$$
\label{bigthm2}
\end{theorem}

\begin{proof} By Proposition \ref{goodgraph}, if $\bl\in \R_k(\A,\K)$, then $\bl\in V_k(\Gamma_\bl,\K)$. The reverse inclusion follows from Corollary \ref{graph-cycles}. We have $\Gamma_\bl\in \NG(\sG,\kk)$ by Theorem~\ref{neighbor} and Corollary~\ref{dimen}. 
\end{proof}

\begin{corollary} $$\R(\A,\K)=\bigcup_{\Gamma\in \NG(\sG,\K)} V(\Gamma,\K).$$
\label{bigthm3}
\end{corollary}

\begin{remark} By \cite{Yuz95}, if $V(\Gamma,\kk)\not = 0$ then $V(\Gamma,\kk)$ is contained in the diagonal subspace $K_0(\Gamma,\kk)=\{\bx\in K \mid \x_{[n]}=0\}.$ In particular, if $V(\Gamma,\kk)\not= 0,$ then $\dim K_0(\Gamma,\kk)\geq 2,$ since $V(\Gamma,\kk)$ contains a resonant pair. If the field \kk\ is replaced by a commutative ring with zero divisors, $V(\Gamma,\kk)$ may not be contained in $K_0(\Gamma,\kk)$ - see \cite{Fa04}.
\label{zerosum}
\end{remark}

A maximal element of the family $\{V(\Gamma,\kk) \ | \ \Gamma \in \NG(\sG,\K)\}$ is officially called the  {\em constituent} of $\R(\A,\K)$ determined by $\Gamma.$ It is often referred to as a ``resonance component" (e.g., local, non-local, global) but this term should be used advisedly in the positive-characteristic setting. If $\K$ is a field of characteristic zero, then the constituents of $\R(\A,\K)$ are precisely the irreducible components, by \cite{LY00}. Indeed this is the case in every example we know of, for any field. But we don't have a general proof that constituents are irreducible in general. By results of the next section, this is a question of irreducibility of intersections of Schubert varieties in special position. 

The next result aids in identifying the constituents among the $V(\Gamma,\K).$

\begin{proposition} Suppose $\Gamma'$ is a subgraph of $\Gamma$ and $K(\Gamma,\K) = K(\Gamma',\K)$. Then $V(\Gamma,\K)\subseteq V(\Gamma',\K)$.
\label{nested}
\end{proposition}

The hypothesis is satisfied, in particular, if $\X_{\Gamma'}(\sG)= \X_\Gamma(\sG).$
If \sG\ has rank two, then \ref{nested} implies $V(\Gamma)$ is maximal if and only if $\Gamma$ consists of $n$ isolated vertices. 

The {\em support} of $V(\Gamma,\K)$  is the set of indices $i$ such that $\l_i\not =0$ for some $\bl\in V(\Gamma,\K)$. The support of $V(\Gamma,\K)$ is precisely the set of non-cone vertices of $\Gamma$, which we will call the support of $\Gamma.$ We say $V(\Gamma,\K)$ is {\em essential} if its support is $[n]$.

Identifying the matroids that support constituents of $\R(\A,\K)$ for \A\ of rank three is one of the more interesting aspects of the theory - see \cite{Yuz04,FY06}.
The two conditions: $\Gamma$ neighborly and $K(\Gamma,\K)$ containing two non-parallel vectors parallel along blocks of $\Gamma,$ are somewhat in opposition. To be neighborly generally requires more edges in $\Gamma,$ but too many edges in $\Gamma$ tends to force vectors in $K(\Gamma,\kk)$ parallel along blocks of $\Gamma$ to be actually parallel. This tension accounts for the rarity of arrangements supporting resonant pairs over fields.

\begin{example} If \A\ is the rank-three braid arrangement of Example \ref{braid1}, then $\R(A,\K)$ has five constituents \cite{F7}. Four of them are supported on the non-trivial rank-two flats $136, 145, 235, 246.$ The neighborly graph for $X=136$  for instance, has edges connecting vertices 2, 4, and 5 to every other vertex; the induced subgraph on vertices 1,3 and 6 has no edges. (See Example \ref{local1}.) 

The fifth constituent is supported on the whole arrangement. For $\Gamma$ having edges 12, 34, and 56 as in Figure~\ref{neighborlyfig}(a), the combinatorial component $V(\Gamma,\K)$ is the subspace of $\K^6$ with basis $\{(1,1,0,0,-1,-1), (0,0,1,1,-1,-1)\}.$ 
\label{global1}
\end{example}

\begin{remark} By \cite{FY06}, if \kk\ has characteristic zero, then $V(\Gamma,\kk)$ is an essential constituent of $\R(\A,\kk)$ if and only if $\Gamma\in \NG(\sG,\kk)$ underlies a {\em multinet} structure on \A. This consists of the partition $\Gamma$ together with an assignment of mutually relatively prime multiplicities to the hyperplanes of \A, so that each block of $\Gamma$ contains the same number of hyperplanes, counting multiplicity, and each codimension-two flat is contained in the same number of hyperplanes from each block of $\Gamma,$ again counting multiplicity. (There is an additional connectivity requirement - see \cite{FY06}.) It is known that the number of blocks can be at most five \cite{FY06,Per06}, and there is only one known example with more than three blocks, namely, the $(4,3)$-net on the Hessian arrangement - see \cite{Yuz04}.
\label{multinet}
\end{remark}
\end{section}

\begin{section}{The line structure of $ V(\Gamma,\K)$}
\label{sect-lines}
Geometers of the early $20^{\rm th}$ used a geometric approach to study skew-symmetric forms (over $\Re$ or $\C$), via null polarities on projective space and projective line complexes \cite{Veblen-Young,Pottman}.  In this section we return to their methods in our more general setting.

Note that $V(\Gamma,\K)$ is a homogeneous variety in $\kk^n.$ Let $\pV(\Gamma,\kk)$ denote its projective image. We will see that $\pV(\Gamma,\kk)$ is the carrier of an algebraic line complex determined by certain projective subspaces associated with $\Gamma$. For background on Grassmannians and Schubert varieties we refer the reader to \cite{Kleim,Ful1,GH}; for line complexes see \cite{Pottman}.

Let $K$ be a \kk-vector space of dimension $k>0,$  and let $\P(K)=\{\K\bx \ | \ \bx\in K-\{0\}\}$ be the projective space of $K$. The standard projective space $\P(\K^k)$ is denoted $\P^{k-1}$. If $D$ is a linear subspace of $K$ we will denote $\P(D)$ by $\bar{D}$ - in particular $\P(K)$ may be written as $\bar{K}.$ For $\bx\in K-\{0\}$ we will write $\bar{\kk\bx}$ as $\bar{\bx}.$ 
If $D$ and $D'$ are subspaces of $K$, let $\bar{D}\ast \bar{D'}=\P(D+D')$. In particular, if $\bar{\bx}$ and $\bar{\bn}$ are distinct points of $\bar{K},$ then $\bar{\bx}\ast \bar{\bn}$ is the line passing through $\bar{\bx}$ and $\bar{\bn}$. 

\subsubsection*{Projective line complexes}  Let $K=\kk^n.$  A line $L=\bar{\bx}\ast\bar{\bn}$ in $\bar{K}$ corresponds to an element of the Grassmanian $\G(2,k)$. The $2 \times 2$ minors $L_{ij}=\x_i\n_j-\x_j\n_i,$ $1\leq i < j \leq d,$ of the matrix $\begin{bmatrix}\bx | \bn\end{bmatrix}$ are called the {\em line coordinates} of $L$. They are determined up to scalar multiple by $L$, independent of the choice of \bx\ and \bn. The Pl\"ucker embedding $L \mapsto [L_{ij} \ : \ 1\leq i<j\leq d]$ identifies $\G(2,k)$ with a $2(k-2)$-dimensional subvariety of $\P^N$, $N=\binom{k}{2} -1.$ The Grassmann-Pl\"ucker relations give a particular set of defining equations for the image $\G(2,k)\subset \P^N$.

A {\em line complex} in $\bar{K}$ is an algebraic subset \L\ of the Grassmannian $\G(2,k)$ under the Pl\"ucker embedding, i.e., a set of lines $\bx \ast\bn$ in $\bar{K}$ given by a system of polynomial equations in the line coordinates $L_{ij}$. The {\em carrier} of a line complex \L\ is the algebraic set $|\L|\subseteq \bar{K}$ of points lying on lines of \L. That is, $|\L|=\bigcup \L$. A {\em variety ruled by lines} is a variety which is the carrier of some line complex.

We are mainly interested in line complexes of the following form. If $D$ is a nontrivial subspace of $K$, set $$\L_D=\{L\in \G(2,k) \ | \ L \cap \bar{D}\not= \emptyset \}.$$ In fact $\L_D$ is a {\em linear line complex}: if $B$ is a matrix whose columns give a basis for $D$, then $L=\bar{\bx}\ast \bar{\bn}\in \L_D$ if and only if all maximal minors of  $\begin{bmatrix}B|\bx|\bn\end{bmatrix}$ vanish. Using the Laplace expansion these minors become linear equations in the $L_{ij}$. 
If $\mathcal D$ is an arrangement of nontrivial subspaces in $K$, let $$\L(\D)=\bigcap_{D\in \D} \L_D.$$ 

\subsubsection*{$\pV(\Gamma,\kk)$ as a ruled variety} Let \A\ be a central arrangement of $n$ hyperplanes. The projective resonance variety $\bR(\A,\kk)$ is clearly ruled by lines: if $\bar{\bl}\in \bR(\A,\kk)$ then $\bl$ belongs to a resonant pair $(\bl,\be),$ and the entire line $\bar{\bl}\ast \bar{\be}$ is contained in  $\bR(\A,\kk).$ Furthermore, $\bl(S)$ is parallel to $\be(S)$ for each block $S$ (= maximal clique) of $\Gamma_{(\bl,\be)},$ which means that there is a point $\bar{\bx}$ on $\bar{\bl}\ast \bar{\be}$ whose restriction $\bx(S)$ to $S$ vanishes. This identifies the underlying line complex for $\bR(\A,\kk)$. For simplicity we study one constituent at a time.

Let $\sG$ be the matroid of \A. Fix a graph $\Gamma\in \NG(\sG,\K)$ and set $K=K(\Gamma,\K)$. Assume $\pV(\Gamma,\K)$ is nonempty.
If $S$ is a block of $\Gamma$, set $$D_S=\{\bx\in K\ | \ \x_i=0 \ \text{for all} \ i\in S\}.$$ The {\em arrangement of directrices} associated with $\Gamma$ is the collection $\D_\Gamma$ of subspaces $D_S,$ where $S$ is a block of $\Gamma$. Note, if $i$ is a cone vertex of $\Gamma$, then every $D\in \D_\Gamma$ is contained in the coordinate hyperplane $\x_i=0$. 

\begin{theorem} Suppose $\bl\in V(\Gamma,\K)$ and $\be\in Z_\Gamma(\bl)-\kk\bl.$ Then, for any block $S$ of $\Gamma$, $\bar{\bl}\ast\bar{\be}$ meets $\bar{D}_S$. Conversely, if $L=\bar{\bl}\ast\bar{\be}$ is a line in $\bar{K}$ which meets $\bar{D}_S$ for every block $S$ of $\Gamma$, then $L\subseteq \pV(\Gamma,\kk)$.
\label{lines}
\end{theorem}

\begin{proof} Since $\bl, \be \in K,$ and $\bl(S)$ is parallel to $\be(S)$, $\bar{\bl}\ast\bar{\be}$ meets $\bar{D}_S$ as observed above. For the converse, suppose \bl\ and \be\ are not parallel, and $\bar{\bl}\ast\bar{\be}$ meets $\bar{D}_S$ for each block $S$ of $\Gamma$. Then for every  block $S$ there exist $a,b\in \K$ such that $a\bl+b\be\in D_S$, consequently $a\bl(S)+b\be(S)=0$. Hence $\be(S)$ is parallel to $\bl(S)$. Since $\bl,\be\in K$ by assumption, this puts $\be$ in $Z_\Gamma(\bl)$. This implies $\bl\ast\be\subseteq \pV(\Gamma).$ 
\end{proof}

\begin{corollary} $\pV(\Gamma)=|\L(\D_\Gamma)|.$
\label{bigthm1}
\end{corollary}

Given a line complex \L\ and $\bar{\bx}\in |\L|,$ the {\em cone} of $\bar{\bx}$ in \L\ is the line complex $\L_{\bar{\bx}}=\{L\in \L \ | \ \bar{\bx}\in L\}$. If $\D=\D_\Gamma$ for some graph $\Gamma$, and $\L=\L(\D),$ then by Theorem~\ref{lines}, $|\L_{\bar{\bx}}|=\P(Z_\Gamma(\bl)).$  In particular  $|\L_{\bar{\bx}}|$ is linear. Define the {\em depth} of $\bar{\bx}$ in \L\ by $\depth(\bar{\bx})=\dim |\L_{\bar{\bx}}|,$ or, equivalently, $\depth(\bar{\bx})=\dim \L_{\bar{\bx}} + 1.$ Recall $V_k(\Gamma,\K)=\{\bl \in V(\Gamma,\K) \ | \ \dim Z_\Gamma(\bl,\K)\geq k+1\}.$ 

\begin{proposition} $\bl \in V_k(\Gamma,\K)$ if and only if $\bar{\bl}$ has depth $k$ in $\L(\D_\Gamma)$.
\end{proposition}

\begin{remark}
Recall the definition of $K_0(\Gamma,\kk)$ from Remark~\ref{zerosum}. If \kk\ has characteristic zero, then by results of \cite{LY00}, $K_0=K_0(\Gamma,\kk)$ is spanned by vectors of the form $u_i-u_j, 1\leq i \leq m,$ where $m$ is the number of blocks of $\Gamma,$ and the $u_i$ are positive integer vectors supported on distinct blocks (see also \cite{FY06}). It follows that each $D\in \D_\Gamma$ is a hyperplane in $K_0,$ and therefore $\pV(\Gamma,\kk)=\bar{K}_0.$ In particular, $\pV(\Gamma,\kk)$ is linear. Moreover, every point of $\bar{K}_0$ has depth $\dim(\bar{K}_0)$. Thus $\bar{K}_0\subseteq \bR_k(\Gamma,\kk)$, where $k=\dim(\bar{K}_0)$ \cite{LY00}. This is false if \kk\ has positive characteristic \cite{MatSuc1}; see Examples~\ref{nonfano} and \ref{deleted1}.
\label{cartan}
\end{remark}

\subsubsection*{Dimension and degree of $\boldsymbol{\pV(\Gamma)}$}
The linear line complexes $\L_D$ are in fact Schubert varieties in $\G(2,k)$. The intersection theory of Schubert varieties can be used to determine the degree of $|\L(\D)|$ if the  $\L_D$ intersect properly, for $D\in \D.$ Our situation is complicated by the fact the arrangement of directrices $\D_\Gamma$ for a neighborly graph $\Gamma$  need not be generic. In all of known examples, the intersection $\L(\D_\Gamma)=\bigcap_{D\in D_\Gamma} \L_D$ is proper, once the ambient space $K$ is replaced by $K_0.$ We have not been able to confirm that this is always the case. 

So we will give a short description of the relevant Schubert calculus, without proof. For this discussion assume \kk\ is algebraically closed. See \cite{Kleim,Harris, Ful1,GH} for background on Schubert varieties, and \cite{Ful2,Ful3} for general intersection theory. 

Fix a complete flag of subspaces $0=K_0\subset K_1\subset \cdots \subset K_k=K.$ For a pair $\sigma=(i_1,i_2)$ of  integers $0\leq i_2\leq i_1\leq k-2$, the associated Schubert variety is the collection $W_\sigma$ of lines $L\in \G(2,k)$ satisfying

\begin{enumerate}
\item $L\cap \P(K_{k-1-i_1})\not = \emptyset,$ and
\item $L\subseteq \P(K_{k-i_2})$.
\end{enumerate}

Schubert varieties associated to different flags are projectively equivalent. By choosing a flag which includes the subspace $D$, we obtain the following.

\begin{proposition} The line complex $\L_D$ is projectively equivalent to $W_{(c(D),0)}$ where $c(D)=\codim(D)-1.$ The codimension of $\L_D$ in $\G(2,k)$ is $c(D).$
\label{cell}
\end{proposition}

Recall the depth of a point $\bar{\bx}\in |\L|$ is the dimension of the carrier of the cone of $\bx$ in $\L.$ If $|\L|$ is irreducible, define the depth of \L\ to be the depth of a generic point on $|\L|$. Then $1\leq \depth(\L)\leq \dim(\L) +1.$

\begin{proposition} If $|\L|$ is irreducible, then  $\dim(|\L|)=\dim(\L) - \depth(\L) +2.$
\label{dimensions}
\end{proposition}


The following observation is adapted from \cite{Pottman}; see also \cite[Example 19.11]{Harris} and \cite[Example 8.3.14]{Ful3}. Let $\D$ be an arrangement of nontrivial subspaces in a $\kk$-vector space $K$, with $\dim(K)=k,$ and $\L=\L(\D).$ Assume without loss of generality that $|\L|$ is not contained in a hyperplane of $K.$ (So, in applying these formulae to $V(\Gamma,\kk),$ we will replace $K$ by a subspace of $K_0.$) Then we can compute the degree of $|\L|$ in the Chow ring of the Grassmannian $\G(2,k).$ 

\begin{theorem}  Suppose $|\L|$ is irreducible, and has dimension $d_0.$ Then the degree of $|\L|$ is given by the intersection formula: $$[W_{(d_0-1,0)}]\cdot [\L]=(\deg |\L|) \,[W_\sigma],$$ where $\sigma=(k-2,k-1-\depth(\L)).$
\label{degree}
\end{theorem}
 
If $\L=\bigcap_{D\in\D} \L_D$ is a proper intersection (meaning each irreducible component has the expected codimension, equal to $\sum_{D\in\D} c(D)$), then $\L$ is rationally equivalent to $\bigcap_{D\in \D} W_{(c(D),0)},$ and the left-hand side of the formula in Theorem \ref{degree} can be computed using the Pieri rule. In this case $d_0-1=\dim(\L)=2(k-2)-\sum_{D\in \D} c(D).$

\begin{corollary} Suppose $\L=\bigcap_{D\in\D} \L_D$ is a proper intersection.
Then the degree of $|\L|$ is given by 
$$[W_{(s,0)}]\cdot \prod_{D\in \D} [W_{(c(D),0)}]=(\deg |\L|) \,[W_\sigma],$$ where $s=2(k-2)-\sum_{D\in \D} c(D)$ and $\sigma=(k-2,k-1-\depth(\L)).$
\label{transdegree}
\end{corollary}

These formulae are currently of limited use, since we have no way to check irreducibility or propriety, or compute depth, except by {\em ad hoc} arguments or computational software (see Example~\ref{hessian}).

\end{section}
\begin{section}{Examples}
\label{sect-examples}

In this section we apply these ideas to several examples. The fine structure revealed in Section \ref{sect-lines} is not apparent in resonance varieties over \C. As we have mentioned in Remark~\ref{cartan}, in this case $V(\Gamma,\C)=\bar{K}_0(\Gamma,\C)$ is linear, and  $\L(\D_\Gamma)$ is the complex of all lines in $\bar{K}_0(\Gamma,\C).$ There is only one known example, supported on an arrangement of rank greater than two, for which $\L(\D_\Gamma)$ consists of more than a single line. Nontrivial line structure emerges over fields of positive characteristic. In particular, we will see that the Hessian arrangement supports a resonance component over $\bar{\Z}_3$ that is a singular irreducible cubic threefold. 

\medskip
With an eye toward the linearity issue, we start with a few elementary observations. These results encompass almost all known examples. Suppose \D\ is a subspace arrangement in a  $\kk$-vector space $K$ of dimension $k>0.$  We define the {\em proper part} $\D_0$ of \D\ by $\D_0=\{D\in \D \ | \ \codim(D)>1\}.$ 

\begin{proposition} \begin{enumerate}
\item $\L(\D)=\L(\D_0)$
\item If $|\D_0|=\emptyset$, then  $|\L(\D)|=\bar{K}$, with every point of depth $(k-1)$.
\item If $\D_0=\{D\}$, then $|\L(\D)_\bl|=\bar{\bl}\ast \bar{D}$ for every $\bl\in K$, and  $|\L(\D)|=\bar{K}.$ Every point of $\bar{K}$ has depth $\dim(D)$.
\item If $\D_0=\{D_1,D_2\}$ with $D_1\not = D_2$, then $|\L(\D)|=\bar{D}_1\ast \bar{D}_2.$ Points of $\bar{D}_1\cap\bar{D}_2$ have depth $\dim (D_1 + D_2)$; all other points have depth $\dim (D_1\cap D_2) +1$.
\end{enumerate}
\end{proposition}

The situation becomes quite simple if $\D$ contains one or more 1-dimensional subspaces. If there is only one such line $\bar{\bl},$ then $\L(\D)=\L_{\bar{\bl}}$ and $|\L(\D)|$ is linear. If there are two or more such lines $\bar{\bl}_1, \ldots, \bar{\bl}_r,$ then $\L(\D)$ is empty or consists of a single line, depending on whether the $\bar{\bl}_i$ are collinear. We will call such elements the {\em poles} of \D. 

By these observations, if the arrangement of directrices $\D_\Gamma$ has fewer than three subspaces of codimension greater than one, or contains an element of dimension one, then the $V(\Gamma,\kk)$ is linear. 

\subsubsection*{Local components in $\mathbf \R(\A,\kk)$} Suppose \A\ is a pencil of $n\geq 3$ lines in the plane. Then $\sG(\A)$ is an $n$-point line.
Let $\Gamma$ be the graph on $[n]$ with no edges. Then $\X_\Gamma(\sG)$ consists of the single (nontrivial) line in $\sG$, and the point-line incidence matrix has rank 1. Then $K=K(\Gamma,\kk)=K_0(\Gamma,\kk)$ has dimension $n-1$. The directrices $D_{\{i\}}$ are all hyperplanes, so $\pV(\Gamma,F)=|\L_{\D_\Gamma}|=\bar{K}$. 

For arbitrary $\sG$, a constituent of $\bR(\A,\kk)$ supported on a flat of rank two is called ``local." For $X\in \X_0(\sG),$ define the graph $\Gamma_X$ by $$\Gamma_X=\{\{i,j\} \ | \ |\{i,j\}\cap X|\leq 1\}.$$ Then every point of $[n]-X$ is a cone vertex, $\X_\Gamma(\sG)=\{X\}$, and $V(\Gamma,\kk)=K(\Gamma,\kk)=K_0(\Gamma,\kk)$ is a linear subspace of dimension $|X|-1$ as above. Suppose $\kk$ is algebraically closed. Then $V(\Gamma,F)$ is  irreducible.  In fact, the proof given in \cite{F7} can be adapted to show that $V(\Gamma,\kk)$ is an irreducible component of $\R(\A,\kk)$ in this case. These are the {\em local components} of $\R(\A,\kk)$. 

\medskip
In the examples below, we give arrangements in terms of their defining polynomials. We label the hyperplanes according to the order of factors in the defining polynomial. We illustrate some of the examples using affine matroid diagrams; the interpretation should be clear, but the reader may consult \cite{Ox} for a detailed explanation (or refer to the explanation of Figure~\ref{neighborlyfig}). We specify graphs by listing their maximal cliques in block notation. We will write vectors over $\Z_2$ as bit strings. Computations were done using {\em Mathematica.}

In almost all known examples, $K_0=K_0(\Gamma,\kk)$ has dimension two. Then $\pV(\Gamma,\kk)=\bar{K}_0$ is a line. For instance, this is the case for the essential constituent for the rank-three braid arrangement, over any field. Referring to Example~\ref{global1}, the arrangement of directrices consists of three (collinear) poles $\bar{\bl}_1=(1,1,0,0,-1,-1), \bar{\bl}_2=(0,0,1,1,-1,-1),$ and $\bar{\bl}_3=(1,1,-1,-1,0,0).$ In preparation for the next example, we point out here  that, in case $\charac(\kk)=2,$ each $\bl_i$ is a sum of characteristic functions of two of the three blocks of $\Gamma$: $\bl_1=110011, \bl_2=001111,$ and $\bl_3=111100.$ 



The next example, which  exhibits higher order, non-local resonance that only appears in characteristic two, was found by D.~Matei and A.~Suciu \cite{MatSuc1}. This example provided the original motivation for the present study. It illustrates that resonance in characteristic two is governed by incidences among submatroids of $\sG$, a phenomenon that has farther-reaching consequences in Example \ref{deleted1} below.

\begin{example} Consider the real arrangement defined by $$Q(x,y,z)=(x+y)(x-y)(x+z)(x-z)(y+z)(y-z)z.$$ Its underlying matroid $\sG$ is the non-Fano plane, with nontrivial lines $$136, \ 145, \ 235, \ 246, \ 347, \ 567.$$ If $\kk$ is a field of characteristic zero, no neighborly graph $\Gamma$ with $\supp(\Gamma)=[n]$  satisfies $\dim K_0(\Gamma,F)\geq 2$. Suppose \kk\ is a field of characteristic two and $\Gamma=127|3|4|5|6,$ as in Figure~\ref{neighborlyfig}(b). Then $\Gamma$ is neighborly, and $\X_\Gamma(\sG)=\X_0(\sG)$. The $6\times 7$  incidence matrix for $\X_\Gamma(\sG)$ has rank 4 over \kk. Hence $K= K(\Gamma,\kk)=K_0(\Gamma,\kk)$ has dimension 3, and $\bar{K}$ is a plane. 
The directrices corresponding to the singleton blocks of $\Gamma$ are lines in this plane, while $D_{127}$ is a pole $\bar{\bl},$ $\bl=0011110.$ Thus $\D_0=\{D_{127}\}$, and $\pV(\Gamma,\kk)=\bar{K}$ as noted earlier. The pole $\bar{\bl}$ has depth two, while all other points on this 2-dimensional component have depth 1. This contrasts with the characteristic-zero situation, in which $\bR_k(\A,\kk)$ consists of the $k$-dimensional components of $\bR(\A,k)$ \cite{MatSuc1}.

\begin{figure}[h]
\begin{center}
\includegraphics[height=4cm]{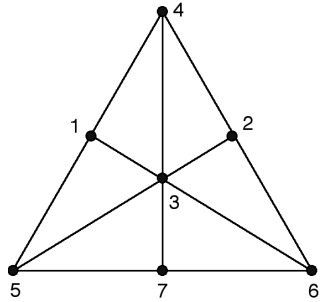}
\end{center}
\caption{The non-Fano plane}
\label{figure-nonfano}
\end{figure}

The diagram of $\sG$  appears in Figure~\ref{figure-nonfano}.  Note that deleting, respectively, points 1, 2, and 7 give submatroids of $\sG$ isomorphic to $\sK_4$, the matroid of rank-three braid arrangement  - see Example~\ref{braid1} and Figure~\ref{neighborlyfig}(a)). Referring to the discussion preceding this example, observe that the special weight $\bl$ is the characteristic function of the intersection of these three submatroids, and in fact is the sum of characteristic functions of two of the three blocks in each of them. For each $i=1,2,7$, there is a resonant pair $(\bl,\be_i)$ supported on $\sG-\{i\}$. Since $\be_1+\be_2+\be_7=0,$ one gets $\dim Z(\bl)=3$. Note that there are no resonant pairs $(\bl,\be)$ with $\be\in (\Z_2)^7$ which are supported on \sG. So $\Gamma=\Gamma_\bl$ does not coincide with any $\Gamma_{(\bl,\be)}$ for $\be\in (\Z_2)^7$ - it is necessary to pass to a field extension (e.g., $\bar{\Z}_2$) to find a generic partner for $\bl.$

The incidence geometry that yields resonance over $\Z_2$ is also reflected in the characteristic varieties $\V_k^1(\A)$: the three components of $\R(\A,\C)$ corresponding to the $\sK_4$ submatroids exponentiate to three 2-tori in $\V^1(\A)\subseteq (\C^*)^7,$ which intersect at the point $(1,1,-1,-1,-1,-1,1)=\exp(2\pi i\bl/2).$ This point is precisely $\V^1_2(\A)$.  This was the first known example of a component of a characteristic variety which does not pass through the identity \cite{CS3}. 
\label{nonfano}
\end{example}

\subsubsection*{The deleted $B_3$ arrangement} In \cite{Suc} A. Suciu introduced the ``deleted $B_3$ arrangement," obtained by deleting one plane from the reflection arrangement of type $B_3$. Suciu showed that the characteristic variety $\V^1(\A)\subseteq (\C^*)^n$ has a one-dimensional component which does not pass through $(1,\ldots, 1)$. Since the components of $\V^1(\A)$ containing  $(1,\ldots, 1)$ are tangent to the components of $\R(\A,\C)$, they all have dimension at least two, by Remarks~\ref{cartan} and \ref{zerosum}. To that point, no arrangements had been found with other than $0$-dimensional components in $\V_1(\A)$ away from $(1,\ldots,1)$. In the next example we see that the same incidence structure that gives rise to this translated component in $\V^1(\A)$ yields components of $\R(\A,\kk)$ with nontrivial intersection, for $\charac(\kk)=2$, in contrast to the characteristic-zero case.
\begin{example} The defining polynomial of the deleted $B_3$ arrangement is given by 
$$Q(x,y,z)=(x+y+z)(x+y-z)(x-y-z)(x-y+z)(x-z)x(x+z)z.$$
Let \sG\ be the underlying matroid, illustrated in Figure \ref{db3diagram}.

\begin{figure}[h]
\begin{center}
\includegraphics[height=4cm]{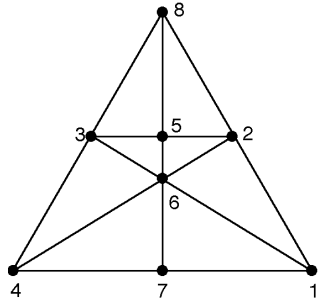}
\end{center}
\caption{The deleted $B_3$ matroid}
\label{db3diagram}
\end{figure}

The deletions $\sG_7=\sG-\{7\}$ and $\sG_5=\sG-\{5\}$ of $\sG$ are copies of the non-Fano plane, and their intersection $\sG_{57}$ is a copy of $\sK_4$. Let $\K$ be an algebraically closed field of characteristic two. With Example \ref{nonfano} in mind, we set $\bl_1=01100101$ and $\bl_2=10010101.$ Then $Z(\bl_1)=V(\Gamma_1,\K),$ where $\Gamma_1=1457|27|37|67|78$ corresponds to the neighborly partition $145|2|3|6|8$ of $\sG_7.$  Similarly, $Z(\bl_2)=V_1(\Gamma_2,\K)$ where $\Gamma_2=2357|15|45|56|58,$ with $\supp(\Gamma_2)=\sG_5.$

Now let $\be=11110000.$ Observe that $\be$ is supported on $\sG_{57}$, and in fact is a sum of characteristic functions of blocks of the neighborly partition $14|23|68$ of $\sG_{57}$. Each of $\bl_1=01100101$ and $\bl_2=10010101$ is a sum of blocks of the same partition. Thus $(\bl_1,\be)$ and $(\bl_2,\be)$ are resonant pairs (Example~\ref{global1}), and so $\be\in Z(\bl_1)\cap Z(\bl_2).$ Indeed,  $Z(\bl_1)\cap Z(\bl_2)=Z(\be).$ In particular $V_1(\Gamma_1,\K)\cap V_1(\Gamma_2,\K)$ is nontrivial. Since $\sG$ itself does not support any neighborly partitions, $V_1(\Gamma_i,\K)$ is indeed an irreducible component of $\R(\A,\K)$ for $i=1,2.$

In Figure \ref{db3picture} is a picture of $\pV(\Gamma_1)\cup \pV(\Gamma_2)$, with the line structure indicated in bold.

\begin{figure}[h]
\begin{center}
\includegraphics[height=6cm]{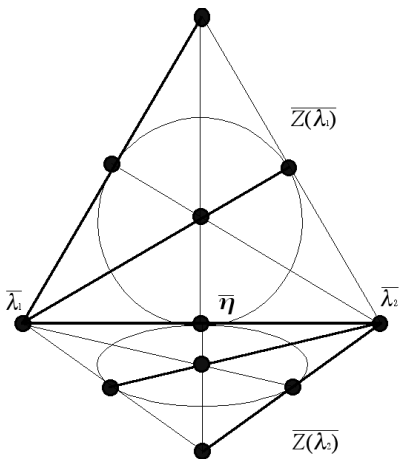}
\caption{A nontrivial intersection}
\label{db3picture}
\end{center}
\end{figure}
Note that $\bl_1\in Z(\bl_2)$ but $Z(\bl_1)\not\subseteq Z(\bl_2).$ That is, there exists $\bx\in \K^7$ such that $a_{\bx} a_{\bl_1}=0=a_{\bl_1} a_{\bl_2},$ but $a_{\bx} a_{\bl_2}\not =0,$ suggesting the possible existence of nontrivial triple Massey products over $\Z_2.$ (Recall that elements of $A_\kk(\A)$ correspond to cohomology classes on $X(\A).$) This phenomenon is expected when $\R_2(\sG,\K)$ is nontrivial and strictly contained in $\R_1(\sG,\K),$ as in Example~\ref{nonfano}. For that example, D.~Matei informs us that there are no nontrivial triple Massey products. We do not know whether the same is true in the present example, where $\R_2(\sG,\kk)=0.$ See \cite{Mat05} for examples of nontrivial $\Z_p$-Massey products on arrangement complements arising from $\Z_p$-resonance varieties.

Referring to \cite{Suc} we find that $\bl_1$ and $\bl_2$ correspond via exponentiation $\bl\mapsto \exp(2\pi i \bl/2)$ to the two points $\rho_1$ and $\rho_2$ that comprise the second characteristic variety $\V^1_2(\A)$, points which lie on the translated component $$C=\{(t,-t^{-1},-t^{-1},t,t^2,-1,t^{-2},-1) \ | \ t\in \C^*\}$$ of $\V_1(\A).$ 
The difference $\bl_1-\bl_2$ is $\be,$ which exponentiates to $\rho_1\rho_2^{-1}.$ Thus $Z(\be)$ exponentiates to the one-dimensional subgroup whose coset by $\rho_i$ is $C$.
The diagram \cite[Figure 6]{Suc} indicates that the same overlapping of non-Fano's and $\sK_4$'s in $\sG$ gives rise to the existence of $\be$ and to that of $C$. 
\label{deleted1}
\end{example}

\subsubsection*{The Hessian arrangement} Finally we present an example having nonlinear components in $\R(\A,\K).$ Let \A\ be the Hessian arrangement in $\C^3,$ corresponding to the set of twelve lines determined by the nine inflection points of a nonsingular cubic in $\P^2(\C)$ \cite[Example 6.30]{OT}.  The underlying matroid is the deletion of one point from $PG(2,3),$ the projective plane over $\Z_3.$
We choose a labelling so that 

$$\X_0(\sG)=\{149\beta, 157\alpha, 168\gamma, 247\gamma,258\beta, 269\alpha, 348\alpha, 359\gamma, 367\beta\},$$

where $H_\alpha, H_\beta,$ and $H_\gamma$ are the last three hyperplanes.
\begin{example} Let $\K=\bar{\Z}_3, $  and let $\Gamma=123| 456 | 789|\alpha\beta\gamma$. Then we have $\X_\Gamma(\sG)=\X_0(\sG).$ The $9 \times 12$ point-line incidence matrix has rank six over $\K$, so $\dim(K(\Gamma,\K))=6.$ The diagonal subspace $K_0(\Gamma,\kk)$ has dimension 5. For each block $S$ of $\Gamma,$ the corresponding directrix has dimension three. Then the projectivized arrangement of directrices  consists of four planes in $\P(\K_0)\cong \P^4.$ (For $\K=\C$ the  directrices are four lines in $\P^2.$)

The placement of these four planes is special: each meets the other three in three collinear points. It follows that the six points of intersection are coplanar, and are the six points of intersection of four lines in general position in that plane. 

A {\em Macaulay2} \cite{M2} computation tells us that the associated ruled variety is an irreducible cubic hypersurface in $\P^4.$ 

This agrees with Schubert calculus formulae from Section~\ref{sect-lines}. 
By Proposition \ref{cell}, a plane $\bar{P}$ in $\P^4$ determines a line complex $\L_P$ equivalent to the Schubert variety $W_{(1,0)}$. The intersection $\L=\L(\D_\Gamma)=\bigcap_{P\in \D_\Gamma} \L_P$ is proper (see below), so $[\L]=[W_{(1,0)}^4]$ in the Chow ring, which in turn equals  $[3W_{(3,1)}+2W_{(2,2)}]$ by the Pieri rule (see \cite{Ful1}). We make the following calculation: 
$$[W_{(2,0)}]\cdot [3W_{(3,1)}+2W_{(2,2)}] =3[W_{(3,3)}],$$
again using the Pieri rule.
Then $\deg |\L|=3$ by Corollary \ref{transdegree}. 

We can show the intersection is proper  by the following {\em ad hoc} argument. The codimension of $\L$ in $G(2,5)$ is at most 4, since that is the codimension of $W_{(1,0)}^4,$  and codimension does not increase under degeneration. One sees without much difficulty that $\L$ has depth one. Then, $\codim \L=\codim |\L|+3$ by Theorem \ref{dimensions}. Since $|\L|$ is easily seen to be a proper subvariety of $\P^4$ we conclude $\codim \L=4,$ as desired. 

Then, using the degree calculation, we can show $|\L|$ is irreducible. Suppose not. Then, being a cubic hypersurface, $|\L|$ contains a $\P^3.$ If three or more of the planes in $\D_\Gamma$ meet this $\P^3$ in lines, necessarily distinct, then $|\L|\cap \P^3$ would be a proper subvariety of $\P^3,$ a contradiction. Thus $\P^3$ must contain at least two of the planes of $\bar{D}_\Gamma.$ Those planes would then meet in a line in $\P^3.$ But the planes in $\bar{\D}_\Gamma$ meet in points. Thus $|\L|$ is irreducible.

H.~Schenck analyzed the scheme structure of cubic threefold $|\L|$ further using {\em Macaulay2} \cite{M2}. The plane containing the six intersection points of the directrices is singular in $|\L|$. These are the points of depth two, i.e, the singular $\P^2$ is a component of $\pV_2(\Gamma,\bar{\Z}_3).$ The quadric in that plane consisting of the four lines containing the intersection points is an embedded component. 
\label{hessian}
\end{example}

It is not hard to see that $\pV(\Gamma,\bar{\Z}_3)$ is the only constituent of $\bR(\A,\bar{\Z}_3)$ that is not linear. Indeed, the only other nonlocal constituents arise from 9-line subarrangements obtained by deleting blocks of $\Gamma$ from \A, and these all give rise to linear constituents.  This is the only known example of a nonlinear resonance constituent over any field. It would be interesting to see a (quadric) resonance constituent that arises from the {\em regulus}, that is, from the line complex determined by three skew lines in $\P^3.$

The variety $\pV=\pV(\Gamma,\bar{\Z}_3)$ from Example ~\ref{hessian} is apparently  related to (the $\bar{\Z}_3$ version of)  more familiar threefolds \cite{Harris,Todd}. It seems to be a fundamental geometric object in the characteristic three universe, arising as it does in such a special way from a classical configuration, and deserves further study. In particular, it is not clear how its existence and structure influences the characteristic variety and fundamental group of the Hessian arrangement.
\end{section}

\begin{ack} This research was helped along by discussions with many people. In the beginning, work by my REU students Cahmlo Olive and Eric Samansky, in the summer of 2000, provided interesting examples to scrutinize. I am grateful to Sergey Yuzvinsky, Hiroaki Terao, Dan Cohen, Alex Suciu, Hal Schenck, and the other participants in two mini-workshops at Oberwolfach (March, 2002 and November, 2003) for their help. Frank Sottile helped me to gain a rudimentary understanding of intersection theory.
\end{ack}


\bigskip
\flushleft{\footnotesize \itshape Department of
Mathematics and Statistics\\ Northern Arizona University \\ Flagstaff, AZ 86011-5717\\
{\tt michael.falk@nau.edu}\\
{\tt http://www.cet.nau.edu/~falk}}

\end{document}